\documentclass[12pt,a4paper]{article}

\usepackage{epsfig}

\usepackage{amsmath}
\usepackage{amssymb}
\usepackage{graphicx}

\newtheorem{propo}{Proposition}[section]
\newtheorem{defi}[propo]{Definition}

\newtheorem{lemma}[propo]{Lemma}
\newtheorem{corol}[propo]{Corollary}

\newtheorem{theo}[propo]{Theorem}

\newtheorem{cor}[propo]{Corollary}

\newtheorem{rem}[propo]{Remark}

\newcommand{\bl}{\begin{lemma}}
\newcommand{\el}{\end{lemma}}

\def\d12{{_{12}}}

 \def\a{\alpha}  
\def\a{\alpha}

\def\Im{{\rm Im}}

\usepackage{indentfirst,latexsym,bm}

\def\Aut{{\rm Aut}}

\def\a{\alpha}

\def\diam{{\rm diam}}

\def\K{{\rm K}}

\def\girth{{\rm girth}}
\def\val{{\rm val}}

\begin{document}
\title{Line graphs and $2$-geodesic transitivity }

\author{Alice Devillers, Wei Jin\thanks{The  second author is  supported by
the Scholarships for International Research Fees (SIRF) at UWA.},
Cai Heng Li and Cheryl E. Praeger\thanks{ This paper  forms part of
Australian Research Council Federation Fellowship FF0776186 held by
the fourth author. The first author was supported by UWA as part of
the Federation Fellowship project during most of the work for this
paper. }\thanks{E-mail addresses: alice.devillers@uwa.edu.au
(A.Devillers),20535692@student.uwa.edu.au
(W.Jin),cai.heng.li@uwa.edu.au (C.H.Li) and
cheryl.praeger@uwa.edu.au (C.E.Praeger).  }
\\{\footnotesize
Centre for the Mathematics of Symmetry and Computation, School of
Mathematics and Statistics,  }
\\{\footnotesize
The University of Western Australia, Crawley, WA 6009, Australia} }

\date{ }

\maketitle

\begin{abstract}
For a graph $\Gamma$, a positive integer $s$ and a subgroup $G\leq
\Aut(\Gamma)$, we prove that $G$ is transitive on the set of
$s$-arcs of $\Gamma$ if and only if $\Gamma$ has girth at least
$2(s-1)$ and $G$ is transitive on the set of $(s-1)$-geodesics of
its line graph. As  applications,  we first prove that the only
non-complete locally cyclic $2$-geodesic transitive graphs are the
complete multipartite graph $K_{3[2]}$ and the icosahedron. Secondly
we classify 2-geodesic transitive graphs of valency 4 and girth 3,
and determine which of them are geodesic transitive.

\end{abstract}

\vspace{2mm}

\section{Introduction}

In this paper, all graphs  are finite simple and undirected. An
\emph{arc} of  a graph  is an ordered vertex pair such that the two
vertices are adjacent. A vertex triple $(u,v,w)$ in a non-complete
graph $\Gamma$ with $v$ adjacent to both $u$ and $w$ is a
\emph{$2$-arc} if $u\neq w$, and a \emph{$2$-geodesic} if the
distance $d_{\Gamma}(u,w)=2$. A graph $\Gamma$ is said to be
\emph{$2$-arc transitive} or \emph{$2$-geodesic transitive} if its
automorphism group $\Aut(\Gamma)$ is transitive on arcs, and on the
2-arcs or 2-geodesics respectively. For connected graphs of girth at
least 4 (where the \emph{girth} is the length of the shortest
cycle), each 2-arc is a 2-geodesic so the sets of 2-arc transitive
graphs and 2-geodesic transitive graphs are the same. However, there
are also connected 2-geodesic transitive graphs of girth 3. It was
shown in \cite[Theorem 1.1]{DJLP-clique} that for such graphs
$\Gamma$, the subgraph $[\Gamma(u)]$ induced on the set $\Gamma(u)$
of vertices adjacent to $u$ is either a connected graph of diameter
2, or is isomorphic to the disjoint union $mK_r$ of $m$ copies of a
complete graph $K_r$ with $m\geq 2,r\geq 2$.

One of the aims of this paper is to characterise 2-geodesic
transitive graphs of girth 3 and valency 4, the smallest valency for
which both possibilities for $[\Gamma(u)]$ arise, namely
$[\Gamma(u)]\cong C_4$ or $2K_2$ for $u\in V(\Gamma)$. This involves
the \emph{line graph} $L(\Sigma)$ of a graph $\Sigma$, namely the
graph whose vertices are the edges of $\Sigma$, with two edges
adjacent in $L(\Sigma)$ if they have a vertex in common.

\begin{theo}\label{smallval-4}
Let $\Gamma$ be a finite connected non-complete graph of girth $3$
and valency $4$. Then $\Gamma$ is $2$-geodesic transitive  if and
only if $\Gamma$ is either $L(K_4)\cong K_{3[2]}$ or  $L(\Sigma)$
for a connected $3$-arc transitive cubic graph $\Sigma$.

Moreover, $\Gamma$ is geodesic transitive  if and only if
$\Gamma=L(\Sigma)$ for a cubic distance transitive graph $\Sigma$,
namely  $\Sigma=K_4$, $K_{3,3}$, the Petersen graph, the Heawood
graph or Tutte's $8$-cage.
\end{theo}

Since there are infinitely many $3$-arc transitive cubic graphs,
there are therefore infinitely many $2$-geodesic transitive graphs
with girth 3 and valency 4. Theorem \ref{smallval-4} provides a
useful method for constructing 2-geodesic transitive graphs of girth
3 and valency 4 which are not geodesic transitive, an example being
the line graph of a triple cover of Tutte's 8-cage constructed in
\cite{Morton-1}. Geodesic transitivity is defined in Section 2. The
line graphs mentioned in the second part of Theorem \ref{smallval-4}
are precisely  the distance transitive graphs of valency 4 and girth
3 given, for example, in \cite[Theorem 7.5.3 (i)]{BCN}. For two
integers $m\geq 3, b\geq 2$, $\K_{m[b]}$ denotes the \emph{complete
multipartite graph} with $m$ parts of size $b$.

One consequence of Theorem \ref{smallval-4} is a classification of
locally cyclic, 2-geodesic transitive graphs in Corollary
\ref{2gt-local-nocycle}: for $[\Gamma(u)]\cong C_n$ is connected and
has diameter 2 only for valencies $n=4$ or 5, and the valency 5,
girth 3, 2-geodesic transitive graphs were classified in
\cite{DJLP}. We note that locally cyclic graphs are important for
studying embeddings of graphs in surfaces, see for example
\cite{Malnic-2,Malnic-1,Malnic-3}. We are grateful to Sandi
Malni$\check{c}$ for  suggesting that we consider $s$-geodesic
transitivity for locally cyclic graphs.

\begin{corol}\label{2gt-local-nocycle}
Let $\Gamma$ be a finite connected, non-complete, locally cyclic
graph. Then $\Gamma$ is  $2$-geodesic transitive   if and only if
$\Gamma$ is  $ K_{3[2]}$ or  the icosahedron.
\end{corol}

Our second aim in the paper is to study further the disconnected
case $[\Gamma(u)]\cong mK_r$ for the smallest value of $m$, namely
$m=2$. Each such graph is isomorphic to the line graph  of some
graph, see \cite[Corollary 1.5]{DJLP-clique}. We investigate
connections between symmetry properties of a connected graph
$\Gamma$ and its line graph $L(\Gamma)$. A key ingredient in this
study is a collection of injective maps $\mathcal{L}_s$, $s\geq 1$,
where $\mathcal{L}_s$ maps the $s$-arcs of $\Gamma$ to certain
$s$-tuples of edges of $\Gamma$ (vertices of $L(\Gamma)$) as defined
in Definition \ref{map-1}. The major properties of $\mathcal{L}_s$
are derived in Theorem \ref{s-arc-tuple-1} and the main consequence
linking the symmetry of $\Gamma$ and $L(\Gamma)$ is given in Theorem
\ref{graph3-lineg2}. This is given in terms of $s$-geodesics and
$s$-arcs (defined in Section 2). The diameter of a graph $\Gamma$ is
denoted by $\diam(\Gamma)$.

\begin{theo}\label{graph3-lineg2}
Let $\Gamma$ be a finite connected regular,  non-complete graph of
girth {\sf g} and valency at least $3$. Let  $G\leq \Aut(\Gamma)$
and let $s$ be a positive integer such that $2\leq s\leq
\diam(L(\Gamma))+1$. Then $G$ is transitive on the set of $s$-arcs
of $\Gamma$ if and only if $s\leq {\sf g}/2+1$ and $G$ is transitive
on the set of $(s-1)$-geodesics of $L(\Gamma)$.
\end{theo}

Note that for the graph $\Gamma$ and the integer $s$ in Theorem
\ref{graph3-lineg2}, there is an additional restriction on  $s$. It
follows from a deep theorem of Richard Weiss in \cite{weiss} that if
$\Gamma$ is $s$-arc transitive, then  $s\leq 7$. This observation
yields  the following corollary.

\begin{cor}\label{graph3-lineg2-cor}
Let $\Gamma$ and   {\sf g} be as in Theorem {\rm \ref{graph3-lineg2}
}. Let $s$ be a positive integer such that $2\leq s\leq
\diam(L(\Gamma))+1$. If $L(\Gamma)$ is $(s-1)$-geodesic transitive,
then either $2\leq s\leq 7$ or $ s>\max \{7,{\sf g}/2+1\}$.
\end{cor}

\section{Preliminaries}

For a graph $\Gamma$, we use $V(\Gamma),E(\Gamma)$, and
$\Aut(\Gamma)$ to denote its \emph{vertex set}, \emph{edge set} and
\emph{ automorphism group}, respectively.  A graph $\Gamma$ is said
to be \emph{vertex transitive}  if the action of $\Aut (\Gamma)$ on
$V(\Gamma)$  is transitive. Vertex transitive graphs are regular in
the sense that $|\Gamma(u)|$ is independent of $u\in V(\Gamma)$, and
$|\Gamma(u)|$ is called the \emph{valency}, denoted by
$\val(\Gamma)$. Regular graphs of valency 3 are called \emph{cubic
graphs}.

A subgraph $X$ of $\Gamma$ is an \emph{induced subgraph} if two
vertices of $X$ are adjacent in $X$ if and only if they are adjacent
in $\Gamma$.  For $U\subseteq V(\Gamma)$, we denote by $[U]$  the
subgraph of $\Gamma$ induced by $U$.

For two vertices $u$ and $v$ in $V(\Gamma)$, a \emph{walk} from $u$
to $v$ is a finite sequence of vertices $(v_{0},v_{1},\ldots,v_{n})$
such that $v_{0}=u$, $v_{n}=v$ and $\{v_{i},v_{i+1}\}\in E(\Gamma)$
for all $i$ with $0\leq i<n$, and $n$ is called the \emph{length} of
the walk. If  $v_{i}\neq v_{j}$ for $0\leq i < j\leq n$, the  walk
is called a \emph{path} from $u$ to $v$.  The smallest integer $n$
such that there is a path of length $n$ from $u$ to $v$ is called
the \emph{distance} from $u$ to $v$ and is denoted by $d_{\Gamma}(u,
v)$. The \emph{diameter} $\diam(\Gamma)$ of a connected graph
$\Gamma$ is the maximum of $d_{\Gamma}(u, v)$ over all $u, v \in
V(\Gamma)$.

\medskip

Let $G \leq \Aut(\Gamma)$ and $s\leq \diam(\Gamma)$. We say that
$\Gamma$ is \emph{$(G, s)$-distance transitive} if, for any $t\leq
s$ and for any two pairs of vertices $(u_1,v_1)$, $(u_2,v_2)$ at
distance $t$,  there exists $g\in G$ such that
$(u_1,v_1)^g=(u_2,v_2)$. If $s$ is equal to the diameter, the graph
is said to be \emph{$G$-distance transitive}.

For a positive integer $s$, an \emph{$s$-arc} of $\Gamma$ is a walk
$(v_0,v_1,\ldots,v_s)$ of length $s$  such that $v_{j-1}\neq
v_{j+1}$ for $1\leq j\leq s-1$. Moreover, a 1-arc is called an arc.
Suppose $G\leq \Aut(\Gamma)$. Then $\Gamma$ is said to be \emph{
$(G,s)$-arc transitive}, if $\Gamma$ contains an $s$-arc, and for
any two $t$-arcs $\alpha$ and $\beta$ where $t\leq s$, there exists
$g\in G$ such that $\alpha^g=\beta$. The study of $(G,s)$-arc
transitive graphs goes back to Tutte's papers \cite{Tutte-1,Tutte-2}
which showed that if $\Gamma$ is a $(G,s)$-arc transitive cubic
graph then $s \leq 5$. About twenty years later, relying on the
classification of finite simple groups, Weiss  \cite{weiss} proved
that there are no $(G,8)$-arc transitive graphs with valency at
least three.  The family of $s$-arc transitive graphs is a central
object in algebraic graph theory, for more work see
\cite{Baddeley-1, IP-1,Praeger-4,Praeger-1993-1,Weiss-1}.

For a graph $\Gamma$ and a positive integer $1\leq s \leq
\diam(\Gamma) $, an \emph{$s$-geodesic} of $\Gamma$ is a walk
$(v_0,v_1,\ldots,v_s)$ of length $s$  such that
$d_{\Gamma}(v_0,v_s)=s$.  It is clear that 1-geodesics are arcs. For
$G\leq \Aut(\Gamma)$, $\Gamma$ is said to be \emph{ $(G,s)$-geodesic
transitive} if, for $1\leq i\leq s$, $G$ is transitive on the set of
$i$-geodesics; further if $s=\diam(\Gamma)$, then $\Gamma$ is said
to be \emph{ $G$-geodesic transitive}. Moreover,  if we do not wish
to specify the group we will say that $\Gamma$ is \emph{
$s$-geodesic transitive} or \emph{ geodesic transitive}
respectively, and similarly for the other properties.  The study of
$s$-geodesic transitive graphs was initiated in \cite{DJLP}, where
the properties of $s$-distance transitivity, $s$-geodesic
transitivity and $s$-arc transitivity were compared.

\medskip

A \emph{maximum clique} of $\Gamma$ is a clique with the largest
possible size. The \emph{clique graph} $C(\Gamma)$ of $\Gamma$ is
the graph with $V(C(\Gamma))=$ $\{$all maximum cliques of
$\Gamma$$\}$, and two vertices are adjacent if and only if they have
at least  one common vertex in $\Gamma$. In particular, if $\Gamma$
has girth at least 4, then $C(\Gamma)$ is the line graph
$L(\Gamma)$.  For example, $L(C_n)\cong C_n$ for $n\geq 3$ where
$C_n$ is the $n$-cycle, and $L(P_{r})\cong P_{r-1}$ for $r\geq 2$
where $P_r$ is the path with length $r$.  The following fact about
line graphs is well-known.

\begin{lemma}{\rm \cite[p.1455]{LB-1} }\label{Linegraph-automorphism-1}
Let $\Gamma$ be a connected graph.   If $\Gamma$ has at least $5$
vertices, then $\Aut(\Gamma)\cong \Aut(L(\Gamma))$.
\end{lemma}

The \emph{subdivision graph} $S(\Gamma)$ of  a graph $\Gamma$ is the
graph with vertex set $V(\Gamma)\cup E(\Gamma)$ and edge set
$\{\{u,e\}|u\in V(\Gamma),e\in E(\Gamma),u\in e\}$. The link between
the diameters of $\Gamma$ and $S(\Gamma)$ was determined in
\cite[Remark 3.1 (b)]{DDP-1}:
$\diam(S(\Gamma))=2\diam(\Gamma)+\delta$ for some $\delta \in
\{0,1,2\}$. Here, based on the above result, we will show the
connection between the diameters of $\Gamma$ and  $L(\Gamma)$ in the
following lemma.

\begin{lemma}\label{diam-lemma-1}
Let   $\Gamma$ be a finite connected graph with $|V(\Gamma)|\geq 2$.
Then $\diam(L(\Gamma))=\diam(\Gamma)+x$ for some $x\in \{-1,0,1\}$.
Moreover, all three values occur, for example, if $\Gamma=K_{3+x}$,
then $\diam(L(\Gamma))=\diam(\Gamma)+x=1+x$ for each $x$.
\end{lemma}
{\bf Proof.} Let  $d=\diam(\Gamma)$, $d_l=\diam(L(\Gamma))$ and
$d_s=\diam(S(\Gamma))$. Let $(x_0,x_2,\ldots,x_{2d_l})$ be a
$d_l$-geodesic of $L(\Gamma)$. Then by definition of $L(\Gamma)$,
each edge intersection $x_{2i}\cap x_{2i+2}$ is a vertex $x_{2i+1}$
of $\Gamma$ and $(x_0,x_1,x_2,\ldots,x_{2d_l})$ is a $2d_l$-path in
$S(\Gamma)$. Suppose that $(x_0,x_1,x_2,\ldots,x_{2d_l})$ is not a
$2d_l$-geodesic of $S(\Gamma)$. Then there is an $r$-geodesic from
$x_0$ to $x_{2d_l}$, say $(y_0,y_1,y_2,\ldots,y_r)$ with $y_0=x_0$
and $y_r=x_{2d_l}$, such that $r<2d_l$. Since both $x_0,x_{2d_l}$
are in $V(L(\Gamma))$, it follows that $r$ is even, and hence
$d_{L(\Gamma)}(x_0,x_{2d_l})=\frac{r}{2}<d_l$ which contradicts that
$(x_0,x_2,\ldots,x_{2d_l})$ is a $d_l$-geodesic of $L(\Gamma)$. Thus
$(x_0,x_1,x_2,\ldots,x_{2d_l})$ is a $2d_l$-geodesic in $S(\Gamma)$.
It follows from  \cite[Remark 3.1 (b)]{DDP-1} that $d_l\leq
d_s/2\leq d+1$.

Now take a $d_s$-geodesic $(x_0,x_1,\ldots,x_{d_s})$ in $S(\Gamma)$.
If $x_0\in E(\Gamma)$, then $(x_0,x_2,x_4,\ldots,x_{2\lfloor d_s/2
\rfloor})$ is a $\lfloor d_s/2 \rfloor$-geodesic in $L(\Gamma)$, so
$d_l\geq \lfloor d_s/2 \rfloor\geq d$. Similarly we see that
$d_l\geq d$ if $x_{d_s} \in E(\Gamma)$. Finally if both
$x_0,x_{d_s}\in V(\Gamma)$, then $d_s$ is even and
$d_{\Gamma}(x_0,x_{d_s})=d_s/2$.  Moreover
$(x_1,x_3,\ldots,x_{d_s-1})$ is a $(\frac{d_s-2}{2})$-geodesic in
$L(\Gamma)$. By  \cite[Remark 3.1 (b)]{DDP-1}, $d_s=2d$, so $d_l\geq
\frac{d_s-2}{2}=d-1$. $\Box$

\section{Line graphs}

Let  $\Gamma$ be a finite connected  graph. For each integer $s\geq
2$, we define a map  from the set of $s$-arcs of $\Gamma$ to the set
of $s$-tuples of $V(L(\Gamma))$.

\begin{defi}\label{map-1}
{\rm Let $\textbf{a}=(v_0,v_1,\ldots,v_s)$ be an $s$-arc of $\Gamma$
where $s\geq 2$, and  for  $0\leq i<s$, let $e_i:=\{v_i,v_{i+1}\}\in
E(\Gamma)$. Define
$\mathcal{L}_s(\textbf{a}):=(e_0,e_1,\ldots,e_{s-1})$.}
\end{defi}

The following theorem  gives some important properties of
$\mathcal{L}_s$.

\begin{theo}\label{s-arc-tuple-1}
Let $s\geq 2$, let $\Gamma$ be a connected graph containing at least
one $s$-arc,  and let $\mathcal{L}_s$ be as in Definition {\rm
\ref{map-1}}. Then the following statements hold.

{\rm (1)}   $\mathcal{L}_s$ is an injective map from the set of
$s$-arcs of $\Gamma$ to the set of $(s-1)$-arcs of $L(\Gamma)$.
Further,  $\mathcal{L}_s$ is a bijection if and only if either
$s=2$, or $s\geq 3$ and $\Gamma \cong C_m$ or $P_n$ for some $m\geq
3,n\geq s$.

{\rm (2)} $\mathcal{L}_s$ maps $s$-geodesics of $\Gamma$ to
$(s-1)$-geodesics of $L(\Gamma)$.

{\rm (3)} If $s\leq \diam(L(\Gamma))+1$, then the image
$\Im(\mathcal{L}_s)$ contains the set $\mathcal{G}_{s-1}$ of all
$(s-1)$-geodesics of $L(\Gamma)$. Moreover,
$\Im(\mathcal{L}_s)=\mathcal{G}_{s-1}$  if and only if
$\girth(\Gamma)\geq 2s-2$.

{\rm (4)} $\mathcal{L}_s$ is $\Aut(\Gamma)$-equivariant, that is,
$\mathcal{L}_s(\textbf{a})^g=\mathcal{L}_s(\textbf{a}^g)$ for all
$g\in \Aut(\Gamma)$ and all $s$-arcs $\textbf{a}$ of $\Gamma$.

\end{theo}
{\bf Proof.} (1) Let   $\textbf{a}=(v_0,v_1,\ldots,v_s)$ be an
$s$-arc of $\Gamma$  and let
$\mathcal{L}_s(\textbf{a}):=(e_0,e_1,\ldots,e_{s-1})$ with the $e_i$
as in Definition \ref{map-1}. Then  each of the $e_i$ lies in
$E(\Gamma)=V(L(\Gamma))$ and $e_k\neq e_{k+1}$ for $0\leq k\leq
s-2$. Further, since $v_j\neq v_{j+1},v_{j+2}$ for $1\leq j\leq
s-2$, we have $e_{j-1}\neq e_{j+1} $. Thus
$\mathcal{L}_s(\textbf{a})$ is an $(s-1)$-arc of $L(\Gamma)$.

Let $\textbf{b}=(u_0,u_1,\ldots,u_s)$ and
$\textbf{c}=(w_0,w_1,\ldots,w_s)$ be  two  $s$-arcs of $\Gamma$.
Then $\mathcal{L}_s(\textbf{b})=(f_0,f_1,\ldots,f_{s-1})$ and
$\mathcal{L}_s(\textbf{c})=(g_0,g_1,\ldots,g_{s-1})$ are two
$(s-1)$-arcs of $L(\Gamma)$, where  $f_i=\{u_i,u_{i+1}\}$ and
$g_i=\{w_i,w_{i+1}\}$ for $0\leq i<s$. Suppose that
$\mathcal{L}_s(\textbf{b})=\mathcal{L}_s(\textbf{c})$. Then
$f_i=g_i$ for each $i\geq 0$, and hence $f_i\cap f_{i+1}=g_i\cap
g_{i+1}$, that is, $u_{i+1}=w_{i+1}$ for each $0\leq i\leq s-2$. So
also $u_0=w_0$ and $u_s=w_s$, and hence $\textbf{b}=\textbf{c}$.
Thus $\mathcal{L}_s$ is injective.

Now we prove the second part. Each arc of $L(\Gamma)$ is of the form
$\textbf{h}=(e,f)$ where $e=\{u_0,u_1\}$ and $f=\{u_1,u_2\}$ are
distinct edges of $\Gamma$. Thus $u_0\neq u_2$, so
$\textbf{k}=(u_0,u_1,u_2)$ is a 2-arc of $\Gamma$ and
$\mathcal{L}_2(\textbf{k})=\textbf{h}$. It follows that
$\mathcal{L}_2$ is onto and hence is a bijection.  If $s\geq 3$ and
$\Gamma \cong C_m$ or $P_n$ for some $m\geq 3,n\geq s$, then
$L(\Gamma)\cong C_m$ or $P_{n-1}$ respectively, and hence for every
$(s-1)$-arc $\textbf{x}$ of $L(\Gamma)$, we can find an $s$-arc
$\textbf{y}$ of $\Gamma$ such that
$\mathcal{L}_s(\textbf{y})=\textbf{x}$, that is, $\mathcal{L}_s$ is
onto. Thus  $\mathcal{L}_s$ is a bijection. Conversely, suppose that
$\mathcal{L}_s$ is onto, and that $s\geq 3$. Assume that some vertex
$u$ of $\Gamma$ has valency greater than 2 and let
$e_1=\{u,v_1\},e_2=\{u,v_2\},e_3=\{u,v_3\}$ be distinct edges. Then
$\textbf{x}=(e_1,e_2,e_3)$ is a 2-arc in $L(\Gamma)$ and there is no
3-arc $\textbf{y}$ of $\Gamma$ such that
$\mathcal{L}_s(\textbf{y})=\textbf{x}$. In general, for $s=3a+b\geq
4$ with $a\geq 1$ and $b\in \{0,1,2\}$, we concatenate $a$ copies of
$\textbf{x}$ to form an $(s-1)$-arc of $L(\Gamma)$: namely
$(\textbf{x}^a)$ if $b=0$; $(\textbf{x}^a,e_1)$ if $b=1$;
$(\textbf{x}^a,e_1,e_2)$ if $b=2$. This $(s-1)$-arc does not lie in
the image of $\mathcal{L}_s$. Thus each vertex of $\Gamma$ has
valency at most 2. If all vertices have valency 2 then $\Gamma\cong
C_m$ for some $m\geq 3$, since $\Gamma$ is connected. So suppose
that some vertex  $u$ of $\Gamma$ has valency 1. Since $\Gamma$ is
connected and each other vertex has valency at most 2, it follows
that $\Gamma\cong P_n$ for some $n\geq s$.

(2) Let $\textbf{a}=(v_0,\ldots,v_{s})$ be an $s$-geodesic  of
$\Gamma$ and let  $\mathcal{L}_s(\textbf{a})=(e_0,\ldots,e_{s-1})$
as above. If $s=2$, then $\mathcal{L}_s(\textbf{a})$ is a 1-arc, and
hence a  $1$-geodesic of $L(\Gamma)$. Suppose that $s\geq 3$ and
$\mathcal{L}_s(\textbf{a})$ is not an $(s-1)$-geodesic. Then
$d_{L(\Gamma)}(e_0,e_{s-1})=r<s-1$ and there exists an $r$-geodesic
$\textbf{r}=(f_0,f_{1},\ldots,f_{r-1},f_r)$ with $f_0=e_0$ and
$f_r=e_{s-1}$. Since $s\geq 3$ and $\textbf{a}$ is an $s$-geodesic,
it follows that $\{v_0,v_{1}\}\cap \{v_{s-1},v_{s}\}=\emptyset$,
that is, $e_0$ and $e_{s-1}$ are not adjacent in $L(\Gamma)$. Thus
$r\geq 2$.  Since $\textbf{r}$ is an $r$-geodesic, it follows that
the consecutive edges $f_{i-1},f_i,f_{i+1}$ do not share a common
vertex for any $1\leq i\leq r-1$,  otherwise $(f_0,\ldots,
f_{i-1},f_{i+1},\ldots, f_r)$ would be a shorter path than
$\textbf{r}$, which is impossible. Hence we have
$f_h=\{u_h,u_{h+1}\}$ for $0\leq h\leq r$. Then
$(u_1,u_2,\ldots,u_r)$ is an $(r-1)$-path in $\Gamma$,
$\{u_1\}=e_0\cap f_1\subseteq \{v_0,v_1\}$ and $\{u_r\}=f_{r-1}\cap
e_{s-1}\subseteq \{v_{s-1},v_s\}$. It follows that
$d_{\Gamma}(v_0,v_s)\leq d_{\Gamma}(u_1,u_r)+2\leq r+1<s$,
contradicting the fact that $\textbf{a}$ is an $s$-geodesic.
Therefore, $\mathcal{L}_s(\textbf{a})$ is an $(s-1)$-geodesic of
$L(\Gamma)$.

(3) Let $2\leq s\leq \diam(L(\Gamma))+1$ and $\mathcal{G}_{s-1}$ be
the set of  all $(s-1)$-geodesics of $L(\Gamma)$. If $s=2$, then by
part (1), each 1-geodesic of $L(\Gamma)$ lies in the image
$\Im(\mathcal{L}_2)$, and hence $\mathcal{G}_{1} \subseteq
\Im(\mathcal{L}_{2})$. Now suppose inductively that $2\leq s\leq
\diam(L(\Gamma))$ and $\mathcal{G}_{s-1} \subseteq
\Im(\mathcal{L}_{s})$. Let $\textbf{e}=(e_0,e_1,e_2,\ldots,e_{s})$
be an $s$-geodesic of $L(\Gamma)$. Then
$\textbf{e}'=(e_0,e_1,e_2,\ldots,e_{s-1})$ is an $(s-1)$-geodesic of
$L(\Gamma)$. Thus there exists an $s$-arc $\textbf{a}$ of $\Gamma$
such that $\mathcal{L}_{s}(\textbf{a})=\textbf{e}'$, say
$\textbf{a}=(v_0,v_1,\ldots,v_{s})$. Since $e_{s}$ is adjacent to
$e_{s-1}=\{v_{s-1},v_s\}$ but not to $e_{s-2}=\{v_{s-2},v_{s-1}\}$
in $L(\Gamma)$, it follows that $e_{s}=\{v_{s},x\}$ where $x\notin
\{v_{s-2},v_{s-1}\}$. Hence $\textbf{b}=(v_0,v_1,\ldots,v_{s},x)$ is
an $(s+1)$-arc of $\Gamma$. Further,
$\mathcal{L}_{s+1}(\textbf{b})=\textbf{e}$. Thus
$\Im(\mathcal{L}_{s+1})$ contains all $s$-geodesics of $L(\Gamma)$,
that is, $\mathcal{G}_{s} \subseteq \Im(\mathcal{L}_{s+1})$. Hence
the first part of (3) is proved by induction.

Now we prove the second part. Suppose that for every $s$-arc
$\textbf{a}$ of $\Gamma$, $\mathcal{L}_{s}(\textbf{a})$ is an
$(s-1)$-geodesic of $L(\Gamma)$. Let ${\sf g}:=\girth(\Gamma)$.  If
$s=2$, as ${\sf g}\geq 3$, then ${\sf g}\geq 2s-2$. Now let $s\geq
3$. Suppose that ${\sf g}\leq 2s-3$. Then $\Gamma$ has a ${\sf
g}$-cycle $\textbf{b}=(u_0,u_1,u_2,\ldots,u_{{\sf g}-1},u_{{\sf
g}})$ with $u_{{\sf g}}=u_0$. It follows that $\mathcal{L}_{{\sf
g}}(\textbf{b})$ forms  a ${\sf g}$-cycle of $L(\Gamma)$. Thus the
sequence $\textbf{b}'=(u_0,u_1,\ldots,u_s)$ (where we take
subscripts modulo ${\sf g}$ if necessary)  is an $s$-arc of $\Gamma$
and $\mathcal{L}_{s}(\textbf{b}')$ involves only the vertices of
$\mathcal{L}_{s}(\textbf{b})$. This implies that
$d_{L(\Gamma)}(e_0,e_{s-1})\leq \frac{{\sf g}}{2}\leq
\frac{2s-3}{2}<s-1$, that is, $\mathcal{L}_{s}(\textbf{b}')$ is not
an $(s-1)$-geodesic, which is a contradiction.  Thus, ${\sf g}\geq
2s-2$.

Conversely, suppose that ${\sf g}\geq 2s-2$. Let
$\textbf{a}:=(v_0,v_1,v_2,\ldots,v_{s})$ be an $s$-arc of $\Gamma$.
Then  $\mathcal{L}_{s}(\textbf{a})=(e_0,e_1,e_2,\ldots,e_{s-1})$ is
an $(s-1)$-arc of $L(\Gamma)$ by part (1).  Let
$\textbf{a}':=(v_0,v_1,v_2,\ldots,v_{s-1})$. Since ${\sf g}\geq
2s-2$, it follows that $\textbf{a}'$ is an $(s-1)$-geodesic, and
hence by (2),
$\mathcal{L}_{s-1}(\textbf{a}')=(e_0,e_1,e_2,\ldots,e_{s-2})$ is an
$(s-2)$-geodesic of $L(\Gamma)$. Thus $z=d_{L(\Gamma)}(e_0,e_{s-1})$
satisfies $s-3\leq z\leq s-1$. There is a $z$-geodesic from $e_0$ to
$e_{s-1}$, say $\textbf{f}=(e_0,f_1,f_2,\ldots,f_{z-1},e_{s-1})$.
Further, by the first part of (3),  there is a $(z+1)$-arc
$\textbf{b}=(u_0,u_1,\ldots,u_z,u_{z+1})$ of $\Gamma$ such that
$\mathcal{L}_{z+1}(\textbf{b})=\textbf{f}$ and we have
$e_0=\{u_0,u_1\}=\{v_0,v_1\}$ and
$e_{s-1}=\{u_z,u_{z+1}\}=\{v_{s-1},v_{s}\}$. There are 4 cases, in
columns 2 and 3 of Table 1: in each case there is a given
nondegenerate closed walk $\textbf{x}$ of length $l(\textbf{x})$ as
in Table 1. Thus $l(\textbf{x})\geq {\sf g}\geq 2s-2$ and in each
case $l(\textbf{x})\leq s+z-1$. It follows that $z\geq s-1$, and
hence $z=s-1$. Thus
$\mathcal{L}_{s}(\textbf{a})=(e_0,e_1,e_2,\ldots,e_{s-1})$ is an
$(s-1)$-geodesic of $L(\Gamma)$.

(4) This property follows from the definition of $\mathcal{L}_{s}$.
$\Box$

\begin{table}[t]\caption{Four cases of $\textbf{x}$ }
\medskip
\centering
\begin{tabular}{|c|c|c|c|c|c|}
\hline Case     &  $(u_0,u_1)$ & $(u_z,u_{z+1})$    &  $\textbf{x}$ & $l(\textbf{x})$ \\
\hline   1    & $(v_0,v_1)$ &  $(v_{s-1},v_{s})$   &
$(v_{s-1},v_{s-2},\ldots,v_2,v_1,u_2,\ldots,u_{z-1},v_{s-1})$ & $s+z-3$  \\
\hline  2    & $(v_0,v_1)$ & $(v_{s},v_{s-1})$   &
$(v_{s},v_{s-1},\ldots,v_2,v_1,u_2,\ldots,u_{z-1},v_s)$ & $s+z-2$ \\
\hline  3    & $(v_1,v_0)$   & $(v_{s-1},v_{s})$ &
$(v_{s-1},v_{s-2},\ldots,v_2,v_1,u_1,u_2,\ldots,u_{z-1},v_{s-1})$ & $s+z-2$ \\
\hline 4    & $(v_1,v_0)$   &  $(v_{s},v_{s-1})$ &
$(v_{s},v_{s-1},\ldots,v_2,v_1,u_1,u_2,\ldots,u_{z-1},v_s)$  &  $s+z-1$ \\
  \hline
\end{tabular}
\end{table}

\begin{rem}\label{smallval4-rem-invariant}
{\rm (i) The map  $\mathcal{L}_s$ is usually not  surjective on the
set of $(s-1)$-arcs of $L(\Gamma)$.  In the proof of Theorem
\ref{s-arc-tuple-1} (1), we constructed an $(s-1)$-arc of
$L(\Gamma)$ not in $\Im(\mathcal{L}_s)$ for any $\Gamma$ with at
least one vertex of valency at least 3.

(ii) Theorem \ref{s-arc-tuple-1} (1) and (3) imply that, for each
$(s-1)$-geodesic $\textbf{e}$ of $L(\Gamma)$, there is a unique
$s$-arc $\textbf{a}$ of $\Gamma$ such that
$\mathcal{L}_s(\textbf{a})=\textbf{e}$. The $s$-arc $\textbf{a}$ is
not always an $s$-geodesic. For example, if $\Gamma$ has girth 3 and
$(v_0,v_1,v_2,v_0)$ is a 3-cycle, then $\textbf{a}=(v_0,v_1,v_2)$ is
not a 2-geodesic but $ \mathcal{L}_2(\textbf{a})$ is the 1-geodesic
$(e_0,e_1)$ where $e_0=\{v_0,v_1\}$ and $e_1=\{v_1,v_2\}$. }

\end{rem}

We are ready to prove Theorem \ref{graph3-lineg2}.

\medskip

\noindent {\bf Proof of Theorem \ref{graph3-lineg2}.} Let  $\Gamma$
be a connected, regular, non-complete graph of girth ${\sf g}$ and
valency at least 3. Then in particular $|V(\Gamma)|\geq 5$, and by
Lemma \ref{Linegraph-automorphism-1}, $\Aut(\Gamma)\cong
\Aut(L(\Gamma))$. Let $G\leq \Aut(\Gamma)$ and let $2\leq s\leq
\diam(L(\Gamma))+1$.

Suppose first that $G$ is transitive on the set of $s$-arcs of
$\Gamma$. Then by Theorem \ref{s-arc-tuple-1} (4), $G$ acts
transitively on $\Im(\mathcal{L}_{s})$. Since $s-1\leq
\diam(L(\Gamma))$, it follows that $L(\Gamma)$ has $(s-1)$-geodesics
and by Theorem \ref{s-arc-tuple-1} (3), $\Im(\mathcal{L}_{s})$
contains all the $(s-1)$-geodesics. Thus $\Im(\mathcal{L}_{s})$ is
the set of $(s-1)$-geodesics of $L(\Gamma)$ and is a $G$-orbit.
Suppose that $s>\frac{{\sf g}}{2}+1$. Let $(a_0,a_1,\ldots,a_{{\sf
g}-1},a_0)$ be a ${\sf g}$-cycle. Then
$(a_0,a_1,\ldots,a_{s-1},a_{s})$ is an $s$-arc. Since the valency of
$\Gamma$ is greater than 2, there exists a vertex $b$ $(\neq a_s)$
adjacent to $a_{s-1}$ such that $(a_0,a_1,\ldots,a_{s-1},b)$ is an
$s$-arc. Since $G$ is transitive on the set of $s$-arcs of $\Gamma$,
there exists $\a \in G$ such that
$(a_0,a_1,\ldots,a_{s-1},a_{s})^\a=(a_0,a_1,\ldots,a_{s-1},b)$, that
is, $a_{s}^\a=b$. As $a_{s}\in \Gamma_{{\sf g}-s}(a_0)$ (the set of
vertices at distance ${\sf g}-s$ from $a_0$) and $a_0^\a=a_0$, we
have $b\in \Gamma_{{\sf g}-s}(a_0)$. Thus there is a $({\sf
g}-s)$-geodesic from $a_0$ to $b$, say $(a_0,b_1,\ldots,b_{{\sf
g}-s-1},b_{{\sf g}-s}=b)$. The walk $(a_0,a_{{\sf g}-1},a_{{\sf
g}-2},\ldots,a_{s-1},b,b_{{\sf g}-s-1},\ldots,b_1,a_0)$ contains a
cycle with length at most $2({\sf g}-(s-1))$. Since $s-1>\frac{{\sf
g}}{2}$, it follows that $2({\sf g}-(s-1))<{\sf g}$ contradicting
that the girth of $\Gamma$ is ${\sf g}$. Thus  $s\leq \frac{{\sf
g}}{2}+1$.

Conversely,  suppose that $s\leq \frac{{\sf g}}{2}+1$ and  $G$ is
transitive on the $(s-1)$-geodesics of $L(\Gamma)$. Then by the last
assertion of  Theorem \ref{s-arc-tuple-1} (3),
$\Im(\mathcal{L}_{s})$ is the set of $(s-1)$-geodesics, and since
$\mathcal{L}_{s}$ is injective, it follows from Theorem
\ref{s-arc-tuple-1} (1) and (4) that $G$ is transitive on the set of
$s$-arcs of $\Gamma$. $\Box$

\medskip

We give a brief proof of  Corollary \ref{graph3-lineg2-cor}.

\medskip

\noindent {\bf Proof of Corollary \ref{graph3-lineg2-cor}.} Suppose
that $ \Gamma, {\sf g}, s$ are as in Theorem \ref{graph3-lineg2} and
that $\Aut(\Gamma)$ is transitive on the $(s-1)$-geodesics of
$L(\Gamma)$ and $s>7$. Then by \cite{weiss}, $\Aut(\Gamma)$ is not
transitive on the $s$-arcs of $\Gamma$ and so by Theorem
\ref{graph3-lineg2}, $s> \frac{{\sf g}}{2}+1$. $\Box$

\section{Two-geodesic transitive graphs that are locally cyclic or locally $2K_2$}

As discussed in Section 1, a graph $\Gamma$ of valency $n$ is
locally cyclic if $[\Gamma(u)]\cong C_n$, and for such a graph to be
2-geodesic transitive (and hence in particular not a complete
graph), $n$ is 4 or 5. Also if $\Gamma$ has valency 4, and $\Gamma$
is 2-geodesic transitive, then $[\Gamma(u)]\cong C_4$ or $2K_2$.
First we treat the case of valency 4, proving Theorem
\ref{smallval-4}. In the proof, we will use the clique graph
$C(\Gamma)$ of $\Gamma$. Recall that $C(\Gamma)$ is the graph with
vertex set of all maximum cliques of $\Gamma$, and two maximum
cliques are adjacent if and only if they have at least  one common
vertex in $\Gamma$.

\bigskip

\noindent {\bf Proof of Theorem \ref{smallval-4}.} Suppose that
$\Gamma$ is a connected non-complete $2$-geodesic transitive graph
of valency 4, and let $A=\Aut(\Gamma)$ and $v\in V(\Gamma)$. Then
$\Gamma$ is arc transitive, and so  $A_v$ is transitive on
$\Gamma(v)$. If $[\Gamma(v)]\cong C_4$, then it is easy to see that
$\Gamma\cong K_{3[2]}$ (or see \cite[p.5]{BCN} or \cite{Cohen-1}).
So we may assume that $[\Gamma(v)]\cong 2K_2$. It follows from
\cite[Theorem 1.2]{DJLP-clique} that $\Gamma$ is isomorphic to the
clique graph $C(\Sigma)$ of a connected graph $\Sigma$ such that,
for each $u\in V(\Sigma)$, the induced subgraph $[\Sigma(u)]\cong
3K_1$, that is to say, $\Sigma$ is a cubic graph of girth at least 4
and $C(\Sigma)$ is in this case the line graph $L(\Sigma)$.
Moreover, \cite[Theorem 1.2]{DJLP-clique} gives that $\Sigma\cong
C(\Gamma)$. In particular, a cubic graph with girth at least 4 has
$|V(\Sigma)|\geq 5$, so by Lemma \ref{Linegraph-automorphism-1},
$A\cong \Aut(\Sigma)$. Now we apply Theorem \ref{graph3-lineg2} to
the graph $\Sigma$ of girth ${\sf g}\geq 4$. Since
$\Gamma=L(\Sigma)$  is 2-geodesic transitive and $2< \frac{{\sf
g}}{2}+1$, it follows from Theorem \ref{graph3-lineg2}  that
$\Sigma$ is $3$-arc transitive. Therefore, $\Gamma$ is the line
graph of a 3-arc transitive  cubic graph.

Conversely, if  $\Gamma \cong K_{3[2]}$, then it is $2$-geodesic
transitive of girth 3. Now suppose that $\Gamma=L(\Sigma)$ where
$\Sigma$ is a 3-arc transitive cubic graph. If $\Sigma$ had girth 3,
then it would be a complete graph, which is not 3-arc transitive.
Hence  $\Sigma$ has girth at least 4. Then $\Sigma$ is locally
$3K_1$, and by \cite[Remark 1.2 (b)]{DJLP-clique},
$C(\Sigma)=L(\Sigma)$ is locally $2K_2$. Thus $L(\Sigma)$ has
valency 4 and girth 3, and hence $L(\Sigma)$ is not $2$-arc
transitive. By Theorem \ref{graph3-lineg2}  applied to $\Sigma$ with
$s=2$,  $L(\Sigma)$ is $2$-geodesic transitive. This proves the
first assertion of Theorem \ref{smallval-4}.

\medskip

Now we suppose  that $\Gamma$ is geodesic transitive. Then $\Gamma$
is distance transitive, and so by Theorems 7.5.2 and 7.5.3 (i) of
\cite{BCN}, $\Gamma$ is one of  the following graphs:
$K_{3[2]}=L(K_4)$, $H(2,3)=L(K_{3,3})$, or the line graph of the
Petersen graph, the Heawood graph or  Tutte's $8$-cage. Further, by
our  argument above, $K_{3[2]}$ is geodesic transitive; by
\cite[Proposition 3.2]{DJLP}, $H(2,3)$ is geodesic transitive. It
remains to consider the last three graphs.

Let $\Sigma$ be the Petersen graph and $\Gamma=L(\Sigma)$. Then
$\Sigma$ is 3-arc transitive, and it follows  from Theorem
\ref{graph3-lineg2}  that $\Gamma$ is 2-geodesic transitive. By
\cite[Theorem 7.5.3 (i)]{BCN}, $\diam(\Gamma)=3$ and $|\Gamma(w)\cap
\Gamma_3(u)|=1$ for each 2-geodesic $(u,v,w)$ of $\Gamma$. Thus
$\Gamma$ is 3-geodesic transitive, and hence is geodesic transitive.

Let $\Sigma_1$ be the Heawood graph and $\Sigma_2$ be  Tutte's
8-cage. Then $\Sigma_1$ is 4-arc transitive and $\Sigma_2$ is 5-arc
transitive, and hence by Theorem \ref{graph3-lineg2}, $L(\Sigma_1)$
is 3-geodesic transitive and $L(\Sigma_2)$ is 4-geodesic transitive.
By \cite[Theorem 7.5.3 (i)]{BCN}, $\diam(L(\Sigma_1))=3$ and
$\diam(L(\Sigma_2))=4$, and hence both $L(\Sigma_1)$ and
$L(\Sigma_2)$ are geodesic transitive. $\Box$

\bigskip

Finally,  we prove Corollary \ref{2gt-local-nocycle}.

\medskip

\noindent {\bf Proof of Corollary \ref{2gt-local-nocycle}.}  Let
$\Gamma$ be a  connected non-complete locally cyclic graph of
valency $n$. Suppose  $\Gamma$ is $2$-geodesic transitive. Then, as
discussed in the introduction, $n=4$ or 5. If $n=4$, then we proved
in  Theorem \ref{smallval-4}, that $\Gamma\cong K_{3[2]}$ and that
$K_{3[2]}$ is indeed 2-geodesic transitive.  If $\val(\Gamma)= 5$,
then by \cite[Theorem 1.2]{DJLP}, $\Gamma$ is isomorphic to the
icosahedron, and this graph is 2-geodesic transitive.  $\Box$


\begin{thebibliography}{hhhh}

\bibitem{LB-1}
L. Babai, Automorphism Groups, Isomorphism, Reconstruction, Handbook
of Combinatorics, the Mit Press, Cambridge, Massachusetts,
Amsterdam-Lausanne-New York, Vol 2, (1995), 1447--1540.

\bibitem{Baddeley-1}
R. W. Baddeley,  Two-arc transitive graphs and twisted wreath
products. {\it J. Algebraic Combin.} {\bf 2} (1993), 215--237.


\bibitem{BCN}
A. E. Brouwer, A. M. Cohen and A. Neumaier, Distance-Regular Graphs,
Springer Ver-lag, Berlin, Heidelberg, New York, (1989).

\bibitem{Cohen-1}
Arjeh M. Cohen,  Local recognition of graphs, buildings, and related
geometries. In Finite Geometries, Buildings, and related Topics
(edited by William M. Kantor, Robert A. Liebler, Stanley E. Payne,
Ernest E. Shult),  Oxford Sci. Publ., New York. {\bf 19} (1990),
85--94.

\bibitem{DDP-1}
A. Daneshkhah, A. Devillers and C. E. Praeger, Symmetry properties
of subdivision graphs, {\it Discrete Math.} (2011),
doi:10.1016/j.disc. 2011.03.031.

\bibitem{DJLP}
A. Devillers, W. Jin, C. H. Li and C. E. Praeger, On distance,
geodesic and arc transitivity of graphs, preprint, 2011, available
at arxiv.org/abs/1110.2235.

\bibitem{DJLP-clique}
A. Devillers, W. Jin, C. H. Li and C. E. Praeger, Clique graphs and
partial linear spaces, in preparation.


\bibitem{IP-1}
A. A. Ivanov, C. E. Praeger,  On finite affine 2-arc transitive
graphs. {\it European J. Combin.} {\bf 14} (1993), 421--444.

\bibitem{Malnic-2}
M. Juvan, A. Malni$\check{c}$ and B. Mohar, Systems of curves on
surfaces, {\it J. Combin. Theory} {\bf B 68} (1996), 7--22.


\bibitem{Malnic-1}
A. Malni$\check{c}$ and B. Mohar, Generating locally cyclic
triangulations of surfaces, {\it J. Combin. Theory} {\bf B 56}
(1992), 147--164.


\bibitem{Malnic-3}
A. Malni$\check{c}$ and R. Nedela, K-Minimal triangulations of
surfaces, {\it Acta Math. Univ. Comenianae} {\bf LXIV 1} (1995),
57--76.

\bibitem{Morton-1}
M. J. Morton,  Classification of 4 and 5-arc transitive cubic graphs
of small girth, {\it J. Austral. Math. Soc.} {\bf A 50} (1991),
138--149.

\bibitem{Praeger-4}
C. E. Praeger,  An O'Nan Scott theorem for finite quasiprimitive
permutation groups and an application to 2-arc transitive graphs,
{\it J. London Math. Soc.} {\bf (2) 47} (1993), 227--239.

\bibitem{Praeger-1993-1}
C. E. Praeger, On a reduction theorem for finite, bipartite, 2-arc
transitive graphs, {\it Australas. J. Combin.} {\bf 7} (1993)
21--36.





\bibitem{Tutte-1}
W. T. Tutte,  A family of  cubical  graphs, {\it Proc. Cambridge
Philos. Soc.} {\bf 43} (1947), 459--474.


\bibitem{Tutte-2}
W. T. Tutte, On the symmetry of cubic graphs, {\it Canad. J. Math.}
{\bf 11} (1959), 621--624.

\bibitem{Weiss-1}
R. Weiss,  s-transitive graphs, {\it  Algebraic methods in graph
theory}, Vol. I, II, (Szeged, 1978), Colloq. Math. Soc. Janos
Bolyai, 25, North-Holland, Amsterdam-New York, (1981), 827--847.



\bibitem{weiss}
R. Weiss,  The non-existence of 8-transitive graphs, {\it
Combinatorica} {\bf 1} (1981), 309--311.




\end{thebibliography}
\end{document}